 

\baselineskip=14pt
\parskip=10pt
\def\Tilde{\char126\relax}
\def\halmos{\hbox{\vrule height0.15cm width0.01cm\vbox{\hrule height
 0.01cm width0.2cm \vskip0.15cm \hrule height 0.01cm width0.2cm}\vrule
 height0.15cm width 0.01cm}}
\font\eightrm=cmr8  
\font\eighttt=cmtt8
\magnification=\magstephalf

\parindent=0pt
\overfullrule=0in

\bf
\centerline
{
WZ THEORY, CHAPTER II
}
\rm
\bigskip
\centerline{ {\it 
Doron ZEILBERGER
}\footnote{$^1$}
{\eightrm  \raggedright
Department of Mathematics, Temple University,
Philadelphia, PA 19122, USA. 
{\eighttt zeilberg@math.temple.edu} \hfill \break
{\eighttt http://www.math.temple.edu/\Tilde zeilberg/   .}
First version: Nov. 10, 1998.  This version: Nov. 16, 1998
(I thank G\"unter Ziegler for several corrections).
Supported in part by the NSF.
A slightly edited version of the transcript of the
(last fifty minutes of a)
videotaped lecture delivered at the MSRI workshop
on Symbolic Computation in Geometry and Analysis, 
Oct. 12-16, 1998. The lecture was given Oct 15, 1998,
9:30-10:30 local time. The videotape is available
from MSRI's website {\eighttt http://www.msri.org}.
A transcript of the first
ten minutes, summarizing my impressions of the
previous talks, was prepared by the MSRI staff, under
the guidance and editorship of Hugo Rossi, and
will appear in a forthcoming MSRI newsletter.
The present article was solicited by Marie-Francoise
ROY, for a volume that she is co-editing. It will
presumably be translated into French and edited.
This original version is exclusive to 
the personal journal of Shalosh B. Ekhad and
Doron Zeilberger and the XXX archives.
} 
}
 
In Erd\H osese, lecturing math is called {\it preaching}.
I am not a preacher, and I am not a proselytizer. I am 
much more than that, I am a prophet. 
In this lecture, I will tell you
how the future of mathematics will look like. 
This I will do by using the WZ theory as an {\it iconic}
example.
 
I am like Moses. Moses was a sinner, like everybody else.
Even Moses was a sinner.  He was punished, along with the 
rest of the desert-generation, not to enter the promised land.
But since Moses was still much more virtuous than his fellow
Israelites, as a compensation, he was allowed, before his death,
to climb to the top of Mt. Nevo and get a glimpse of the
promised land. I am also doomed to die, along with you, before
the promised land of computer-math will come, but I have a glimpse
of the future that I wish to share with you, fellow mortals.
 
One of the messages I got from this conference was that
`Real' Analysis is really algebra! As I am sure you know,
Brouwer, Weyl and Bishop
have already pointed out that `real' analysis
is a misnomer. It is not at all real, like `complex' analysis
is not at all complex (it is simpler than real analysis).
Real analysis rests on the dubious notion of {\it real number},
that uses the actual infinity and is not constructive.
So traditional real analysis is purely fictional. But 
a large part of analysis is really algebra in disguise.
For example the differentiation operator $D$ can be defined
completely algebraically by Leibnitz's rule $D(xf)=xDf+f$.
When we teach calculus, we really teach algebra, symbolic
manipulations of functions and derivatives.
It is that part of analysis that is amenable to computer algebra,
and that is meaningful. Two celebrated examples are
differential Galois theory, one of whose gurus is Michael
Singer, and the theory of D-modules and holonomic systems,
that is behind the scenes of WZ-theory.
 
Recall one of Dave Bayer's many aphorisms:
 
ALGEBRA WITH PERSONALITY=COMBINATORICS.
 
But the converse is also true, Combinatorics is really algebra,
non-commutative algebra, since you can always encode combinatorial
objects into a computer, and everything a computer can handle
is ultimately algebra.
 
Another corollary is
 
MATH=SYMBOLIC COMPUTATION
 
so ultimately, everything is computer algebra.
This has always been true, except that until recently,
the computer was God-made, that tiny CPU between our shoulders,
some are better than others of course, but more or less they
are all equivalent.
 
Another dichotomy is Symbolic computation vs. Numeric computation.
Symbols are profound, while numbers are trivial. But 
WHAT ARE NUMBERS? they are symbols!, conversely, every symbol
is a number, e.g. its ASCII or uni-code representation. So
numbers are symbols and symbols are numbers, and of course, ultimately
everything is just 0s and 1s, and it all depends on your perspective.
 
So the point is, Don't Be a Snob! and don't look down on
what other people are doing. I just bought a book by Kary Mullis,
the `cool' 1993 Chemistry Noble prize winner, and there he said that
chemists look down on biologists, 
physicists look down on chemists,
mathematicians look down on physicists, and philosophers 
{\it used to} look down on mathematicians,
until they realized that they can do nothing. Another
book that I have just bought, by Don Knuth, puts it more mildly.
Knuth believes that mathematicians have different thought-processes
than physicists, who in turn, think differently than chemists,
and computer scientists think yet differently: algorithmically.
 
But even within math we have this snob-hierarchy.
The theoretical crowd looks down on the computational crowd,
pure looks down on applied etc. Nowadays, the opposite is also
true, so we have mutual contempt. Let's hope that this will give
way to mutual respect and understanding!
 
{\bf WZ Theory, Chapter 0}
 
The trend in mathematics is starting to go from
{\it computer-assisted} to {\it computer-generated}.Proof 
First there are
 
{\it (i) Computer-assisted conjectures }
 
This is nothing new, and goes back to Pythagoras,
Archimedes, Euler, Gauss, Riemann and all the other
giants, who did extensive experimentation to find
conjectures ( e.g. the prime-number theorem and
the Riemann Hypothesis). But, of course, the computer
was God-made. With Man-made computers, people are able
to discover many new amazing facts.
 
One of my favorite recent conjectures, that has been driving me
up the wall for the last five years is
\eject
{\bf Mark Haiman's Macaulay-Generated Conjecture}
 
Let K be a field (of char. $0$), and let $I$ be the
ideal in $K[x_1, \dots, x_n; y_1, \dots, y_n]$ generated
by the {\it non-constant} diagonally symmetric polynomials
(i.e. polynomials $p(x,y)$ such that for all $1 \leq i <n$,
$$p(x_1, \dots, x_{i+1},x_i, \dots, x_n;
   y_1, \dots, y_{i+1},y_i, \dots, y_n)= 
$$
$$p(x_1, \dots, x_i,x_{i+1}, \dots, x_n;
   y_1, \dots, y_i,y_{i+1}, \dots, y_n) \quad .
$$
Then
$$
dim(K[x_1, \dots, x_n; y_1, \dots, y_n]/I)=(n+1)^{n-1} \quad .
$$
 
Macaulay, at the hands of Mark Haiman, proved this for
$n=1, \dots, 7$, but no one, as yet was able to prove it
for {\it all} $n$. Who knows? perhaps the 2010 version of
Macaulay could? By the way, I was thrilled to meet the
developers of Macaulay, Dave Bayer and Mike Stillman, and the
guru behind it, David Eisenbud, since Macaulay is one of my
favorite systems.
 
Next we have
 
{\it (ii) Computer-Generated Conjectures}
 
The mode of doing computer-assisted conjectures, is {\it on-line}.
You sit next to the terminal, drink your cup of coffee
(being careful not to spill on the keyboard), and look at what
the screen spits out. Then you try this and try that, and
have a {\it dialog}, until you arrive at a meaningful conjecture.
But suppose that you don't like to drink coffee, and you  would like
the computer to do all the conjecturing for you. That's possible
too. There exist powerful software that automatically finds
conjectures, but without proving them. Two such excellent systems
are the Maple packages {\it gfun} and {\it mgfun}, developed
by Bruno Salvy and Paul Zimmermann at INRIA, France.
 
We also have
 
{\it (iii) Computer-Assisted Proofs}
 
Many proofs nowadays are computer-assisted, but in most of
them computers are not mentioned. If it is at all possible to hide
the role of the computer, it is not mentioned. Since computers
can't yet sue you, it is safe to cheat and pretend that it
was all pure human insight. Of course, in some cases it is
impossible to hide the role of the computer, like for example,
in Lanford et. al.'s proof of the Feigenbaum conjecture. Since
there, a computer was also needed to verify the proof. Other
obvious proofs where the role of the computer had to be
openly admitted are the Appel-Haken
proof of the 4-Color Theorem and Tom Hales's proof of the Kepler
conjecture. But, if it is at all possible, most people
cheat, and do not mention the computer. (Of course, this
continues a long tradition of mathematicians hiding their
motivations and heuristics, and only presenting the
end product, the proof.)
 
Now the next revolution (Ta Ta) is
 
{\it (iv) Computer-Generated Proofs}
 
Although it may have happened before, the first 
full-fledged computer-generated proofs started with
WZ theory. Here a computer writes, all by itself, the
full publishable article. For example, my beloved computer,
Shalosh B. Ekhad, already has 18 articles, 16 of which are
mentioned in Math Reviews. It could have gotten tenure at any
top-twenty university, maybe not top-ten, but definitely top-twenty,
had it been a human. But there are still some prejudices left,
and at present it has no academic position. Is any one
here hiring by any chance?
 
But how does Shalosh write a paper? It has an algorithm that
generates the paper. It still needs a human to design and
program an algorithm that would generate papers. So we
still need a human behind the scenes, and that human was
me, in this case. But this is only a very intermediate
stage. Soon we would have
 
{\it (v) Computer-Generated Algorithms [but still using human-
generated theories]}
 
to make me superfluous. After that would come
 
{\it (vi) Computer-Generated Theories [but using human-generated,
more general, theories]},
 
but, after awhile, we would finally have
 
{\it (vii) Computer-Generated Theories [using computer-generated
 theories]}.
 
When this time comes,
we humans would have to retire, as principal mathematicians, and
let computers do all the serious math. You might say, but we
still need humans to develop meta-theories, and meta-meta-theories...
But this is an illusion. As any meta-mathematician can tell you,
the meta-operator is idempotent, i.e.:
 
META META=META .
 
So once computers will get that far, they would be self-sufficient,
and do all the `serious' math.
 
But this does not mean that mathematics will be completely
abandoned by humans. We still lift weights, run, jump, and so
on, even though machines can lift much heavier weights, cars
can move much faster, and airplanes can jump much higher.
Human mathematics is going to become a sport, and we'll do it
as recreation. In fact, it will be a competitive sport, and
its stars would become rich celebrities. This is not necessarily
a bad thing. Look how much a baseball star is making today,
and how much, say, David Eisenbud is making. Even if you take
the log, even the loglog, it is significantly more.
 
But this is the future. We are still very human-centrists.
Even that paradigm of political correctness, the TV program
Sesame Street, that advocates that everyone is equal,
regardless of sex, race, color, creed, religion, and even
sexual orientation, once said:
 
``Computers can't think for themselves, they only do what
we tell them.''
 
Other quotations from well-meaning but ignorant HCPs
(Human Chauvinist Pigs) are:
 
``Computers are {\it just} machines for doing faster what we
already know how to do slower'' (--- Gian-Carlo Rota).
 
``The computer is {\it just} a pencil with power-stirring''
(--- George Andrews),
 
and finally,
 
``I never use a computer'' (--- Andrew Wiles).
 
Unlike the previous quotations, Wiles's is entirely factual,
and is not a value judgment. It was uttered at the
Nova program on FLT, that later became the book
`Fermat's Enigma'. I completely trust Wiles that he 
did not `cheat' and use a computer. But the commentary
that followed, by the producer and author Simon Singh,
expressed some already obsolete humo-centrist
conventional wisdom. Singh said:
 
``Computers can only prove the Taniyama-Shimura conjecture
for {\it each specific} elliptic curve, but it takes
a {\it human genius} to prove it for all of them, using
{\it logical deduction}''.
 
Well here is a proof of the irrationality of $\sqrt{2}$ that
could have easily been found by computer.
 
{\bf Theorem}: There are no positive integers $A$ and $B$ such that
$A^2-2B^2=0$.
 
{\bf Proof}: Let $a=2B-A, b=A-B$. Since
$a^2-2b^2=-(A^2-2B^2)$ (check!), any solution $(A,B)$ entails
a smaller solution $(a,b)$. Since $(1,1)$ is not a solution, we have
a contradiction.
 
This trivial proof was found by hand, but it is very conceivable that
one day Shalosh B. Ekhad or one of its siblings will find
a transformation $(x,y,z,n) \rightarrow (x',y',z',n')$ that
preserves Fermaticity, and that shrinks according to some
norm (of course, it would only be valid for $n' \geq 3$).
This putative transformation, only a computer would be able to find, 
and only a computer would be able to verify, 
but modulo routine checking, would be much more elegant
than Wiles's baroque proof. In fact, such a proof was
announced by Shalosh B. Ekhad, on April 1, 1995
(see its article, `Proof of the Riemann Hypothesis
and some other hitherto undemonstrated theorems', available
from its website), but unfortunately its hard disk crashed...,
so finding an elementary proof of FLT is open again.

So here is the fallacy. Conventional wisdom asserts that 
empirical science,
and computer experimentation, use {\it induction} while math uses
{\it deduction}. But induction and deduction got married
a long time ago, via the principle of {\it complete mathematical
induction}. To prove $A_n$ for all $n$, all you have to do is
prove the two statements $A_0$ and $A_n \Rightarrow A_{n+1}$,
both of which can (at present sometimes, in the future always) be
proved by computer.
 
At this rate, I don't know whether I'll ever make it to
Ch. II, but now we can at least have a quick review of
 
{\bf WZ Theory, Chapter I}
 
Here is a deep
 
{\bf Theorem}:
$$
\sum_{k=0}^{n} {{n} \choose {k}} =2^n \quad .
$$
 
According to Simon Singh, a computer can verify this
for $n=0$, $n=1$, $n=2, \dots $ , even for $n=10000000$, but
it takes a human (albeit not necessarily as smart as
Andrew Wiles) to prove it for every $n$. Well, Simon,
I have news for you. My beloved computer, Shalsoh B. Ekhad,
that is all chips and silica, has a one-line proof. To wit:
 
{\bf Proof:} $k/2(n-k-1)$ \halmos
 
This is an almost trivial, {\it iconic} (to use
Dave Bayer's beautiful language) example of so-called
WZ theory. By the way, the `W' in WZ stands for `Wilf'.
I am not sure about the `Z': perhaps Zeno?, Zorn?, Zelmanov?,
don't really know.
 
Of course, we have to explain the {\it WZ methodology} once and
for all. But once we do that, we would be able to understand
and appreciate Shalosh's one-line proofs. 
 
{\bf Explanation for the non-WZist}
 
Whenever we want to prove an identity of the format
$$
\sum_{k} PRETTY(n,k) = NICE(n) \quad,
$$
the first, crucial, step, is to divide by the right hand side,
getting:
$$
\sum_{k} PRETTY(n,k)/NICE(n) =  1 \quad .
$$
Now, PRETTY over NICE is not just `Pretty Nice', it is
super-nice. But since I don't like superlatives, let's
rename the summand $NICE(n,k)$.
 
{\bf Fundamental Theorem of WZ theory:} `Whenever' you have to
prove an identity of the form
$$
\sum_{k} NICE(n,k)=1, \quad, for\,\,\, all\,\,\,n \geq 0 \quad ,
$$
there exists another nice function, $NICE'$, such that (Ta Ta!)
$$
NICE(n+1,k)-NICE(n,k)=NICE'(n,k+1)-NICE'(n,k) \quad .
\eqno(WZ)
$$
If you had to preserve the
top-ten most influential formulas of the 20th century,
$E=mc^2$ would be amongst them,  and so would $(WZ)$!
 
{\bf Why Does the Mere Existence of NICE' Prove The Identity?}
 
If you are skeptical, you can check
at any instance, either by hand or using
your own PC, 
that the $NICE'$ outputted by EKHAD satisfies, with
the inputted $NICE$,
the purely routine identity $(WZ)$. Then, defining
$$
a(n):=\sum_{k} NICE(n,k) \quad,
$$
we have
$$
a(n+1)-a(n)=\sum_{k} \,\,[\,NICE(n+1,k)-NICE(n,k)\,]\,=\,
\sum_{k}\,\, [\,NICE'(n,k+1)-NICE'(n,k)\,]\,=\,0 \quad,
$$
since $a(0)=1$ (check!), if follows that $a(n)=1$ for all $n \geq 0$.
 
{\bf How did the Computer Find NICE'?}
 
By an algorithm! The referee checked it, even Knuth, the
great guru of computer science, checked it carefully,
before he covered it in the classic `Concrete Math' by
Graham, Knuth, and Patashnik.
If you don't trust Knuth, if you think that he is a flake,
go ahead, and check it. But I advise you strongly to believe
the algorithm, since you have probably better things to do.
So just push the button, and wait for the output.
 
There is one more important fact. The quotient
$NICE'(n,k)/NICE(n,k)$ is {\it always} a humdrum,
pedestrian, rational function $R(n,k)$, that we dubbed
{\it certificate}. From $R(n,k)$ 
you can reconstruct $NICE'(n,k):=NICE(n,k)R(n,k)$, 
and if skeptical verify $(WZ)$. But the whole proof
is {\it encapsulated} by $R(n,k)$, and from now on, in order
to prove any hypergeometric identity, all the computer has to do is
output the certificate $R(n,k)$, and that's the {\it the whole proof}.
 
So it is always a {\it one-line proof} which Shalosh B. Ekhad can
find just like that, usually in a few nano-seconds. Sometimes in
a few hours, but usually very fast.
 
If you think that the binomial theorem is not very impressive, here
is another example. Once upon a time there was a poor clerk in India
who was very good in math, but did not pass the English entrance
examination to college, so he had to get a job. Luckily,
his boss was a nice guy, and let him spend lots of time
discovering lots of formulas. Then one day he gathered enough
courage to write a letter to the leading analyst of the day.
Among his thousands of formulas, he had to choose the
most impressive. One of the dozen-or-so formulas that
made it (Ramanujan's letter to Hardy) was:
$$
{{2} \over {\pi}} =
\sum_{k=0}^{\infty} (-1)^k (4k+1){{ (1/2)_k^3} \over
{k!^3}} .
$$
This, by itself, is not-yet-WZable, but let's cheat and
prove the more general identity, that is completely shaloshable.
 
{\bf Theorem:}
$$
{{\Gamma(3/2+n)} \over {\Gamma(3/2) \Gamma ( n+1) }} =
\sum_{k=0}^{\infty} (-1)^k (4k+1) {{ (1/2)_k^2 (-n)_k
} \over
{k!^2 (3/2 + n)_k }}.
$$
{\bf Proof:} ${{-2k^2} \over {(n-k+1)(4k+1)}}$. \halmos
 
To deduce the original identity, we simply "plug" in
$n= -1/2$, which is legitimate in view of Carlson's theorem.
 
So here is a completely {\it computer-generated} proof of
an hitherto non-trivial identity.
 
But here is:
 
{\bf An Even More `Impressive' Example: High-Brow Meets Low-Brow}
 
Ken Ono, one of the most promising young number-theorists
today, working with the even younger, and just as promising,
Scott Ahlgren, tried to prove the following conjecture
of Frits Beukers.
Let
$$ 
A(n):=\sum_{k=0}^{n}{n\choose k}^2{n+k\choose k}^2 \quad ,
$$
be the famous Ap\'ery numbers, and define integers $a(n)$ by
$$ 
\sum_{n=1}^{\infty}a(n)q^n:=q\prod_{n=1}^{\infty}
   (1-q^{2n})^4(1-q^{4n})^4=q-4q^3-2q^5+24q^7- \cdots
$$
Beukers conjectured that
if $p$ is an odd prime, then
$$
  A((p-1)/2 ) \equiv a(p) \pmod{p^2} \quad.
$$
Ken and Scott used many high-powered tools on modular forms,
and were able to reduce it to a lowly binomial-coefficient
identity:
$$
  \sum_{k=1}^{n}k{n\choose k}^2{n+k \choose k}^2
  \left \{ {{1}\over {2k}}+\sum_{i=1}^{n+k}{{1}\over {i}}+
 \sum_{i=1}^{n-k}{{1}\over {i}}-2\sum_{i=1}^{k}{{1}\over {i}}\right \}=0.
$$
And they got stuck. (Often when you are working on a hard problem,
after doing fancy and high-brow stuff, you sooner or later
hit rock, and none of the abstract nonsernse stuff can bail you out.
And that rock is either combinatorics or computer algebra.)
So they got stuck, and sent it to me, offering to collaborate.
I would have been
stuck too..., but I  did not even try, and gave it right away to
my beloved servant Shalosh B. Ekhad, who proved it in a few seconds.
The one-page article,
authored by Ahlgren, Ekhad, Ono, and myself,
has recently appeared in the Electronic Journal of Combinatorics.
 
Now, let's have a short interlude, and discuss the
\eject
{\bf Infra-Structure}
 
Myself, I was brought up 
in the desert-generation, and was brainwashed that everything
is subsidiary to theory. In fact, to be honest, while, by hindsight,
most of WZ theory could have been developed completely elementarily,
it did originate in the context of $D-$modules, holonomic
systems, and elimination in the Weyl algebra
$K[x_1, \dots , x_n; D_1, \dots, D_n]$.
 
{\bf Growing Research Community}
 
In addition to the activity in Philadelphia (both Penn and
Temple), there is a very active center at RISC-Linz, Austria,
under the leadership of Peter Paule. His group consists
of the brilliant students Axel Riese, Markus Schorn,
Kurt Wegscheider, and Burkhard Zimmermann. In INRIA, France,
there is Frederic Chyzak, who beautifully implemented
the holonomic paradigm. In Japan, Nobuki Takayama developed
the marvelous package KAN. Let me take this opportunity and
advertise a
 
{\bf COMING ATTRACTION}
 
Peter Paule's talk this Friday, at 2:00PM. He will talk about
`Fine-Tuning the WZ-engine'.
 
Still in commercial mode, let's take a
 
{\bf COMMERCIAL BREAK}
 
BUY A=B by Petkovsek, W, and Z!
 
If you don't mind paying the list price, \$39.00 you
can order it from: {\tt http://www.amazon.com}.
If that's too much, you can get a new copy for only
\$27.50 at  {\tt http://www.barnesandnoble.com}.
 
If you can't afford even that, then you may be able to find
it, for \$15.00 at the used-books website:
{\tt http://www.bibliofind.com} . I don't know if
they had more than one copy, because if they did,
I have already bought it myself, but perhaps they have
more.
 
Here are excerpts from sample reviews:
 
`Masterpiece'--Peter Paule (Math Reviews)
 
`Great example of mathematical exposition'--
Noam Zeilberger ({\tt http://www.home.com/\Tilde ironblaze/}).
 
The outlier review was:
 
`Sloppily written in parts and lacking implementation details'-
Wolfram Koepf.
 
Luckily, Koepf has remedied the shortcomings of A=B by writing
his own book: `Hypergeometric Summation' that is highly
recommended in case you are willing to buy two books on the
subject.
 
Before starting chapter $1 {{1}\over {2}}$ let me just mention that
whatever I said about sums is all true for integral identities
$$
\int PRETTY(x,y) dy= NICE(x) \quad .
$$
The theory was worked out by Gert Almkvist and myself, published
in 1990 in J. Symbolic Computation. Michael Singer was
the editor, and made it even better.
 
Now we are ready for
 
{\bf WZ Theory, Chapter ${\bf 1 {{1}\over {2}}}$}\footnote{$^2$}
{\eightrm W and Z, Inv. Math. 108(1992), 575-633. Usually I am
too timid to submit my concrete stuff to the fancy and
snooty Inventiones, and I am sure that if Herb and I
had submitted it in the canonical way, to the editor-in-chief
Remmert, we would have gotten the usual rejection-slip form
letter. It so happened, that during the previous semester, one of the
other editors, Marcel Berger, visited Penn, and befriended
Wilf, who even took him on his airplane. So he owed him a favor...}

Not only is the binomial theorem
$$
(x+y)^n=\sum_{k} {{n} \choose {k}} x^k y^{n-k} \quad,
$$
fully shaloshable, but so is the trinomial theorem
$$
(x+y+z)^n=\sum_{k_1,k_2} {{n!} \over {k_1!k_2!(n-k_1-k_2)!}} 
x^{k_1} y^{k_2} z^{n-k_1-k_2} \quad,
$$
and, in principle, the centonial and even zillionomial theorem.
These, and of course much more complicated ones are now
routinely provable thanks to the
 
{\bf Fundamental Theorem of Multi-WZ Theory (Discrete Version)}
 
`Whenever' we want to prove an identity of the form:
$$
\sum_{k_1}\dots  \sum_{k_r} NICE(n; k_1, \dots , k_r) \equiv 1 \quad,
$$
There exist $NICE'_1, NICE'_2, \dots, NICE'_r$ such that
$$
\Delta_n NICE = \sum_{i=1}^{r} \Delta_{k_i} NICE'_i \quad
$$
Furthermore, $R_i:=NICE'_i/NICE$, ($i=1, \dots, r$),
are rational functions of $(n; k_1, \dots , k_r)$.
 
We also have:
{\bf Fundamental Theorem of Multi-WZ Theory (Continuous Version)}
 
`Whenever' We Want To Prove An Identity Of The Form:
$$
\int \dots \int NICE(n; x_1, \dots, x_r) \,\, dx_1 \dots dx_r \equiv 1 \quad,
$$
there exist $NICE'_i(n; x_1, \dots, x_r)$ ($i=1, \dots, r$), such that
$$
\Delta_n NICE = \sum_{i=1}^{r} {{\partial} \over {\partial x_i}} NICE'_i
\quad .
$$
Furthermore, $R_i:=NICE'_i/NICE$ are rational functions of
$(n; x_1, \dots, x_r)$.
 
Note that the mere existence of  $NICE'_i$
(or equivalently the $R_i$), together with
the trivially verifiable case $n=0$, implies the identity, 
since,
$$
\Delta_n \int \dots \int NICE \,=\,
\int \dots \int \Delta_n NICE \,=\,
\sum_{i=1}^{r} \int \dots \int {{\partial} \over {\partial x_i}} \{ NICE'_i\}
=0 \quad ,
$$
(because the $NICE'_i$ are of compact support).
 
Let me recall two `famous' examples. They are not really famous, but
they are due to famous people.
The first one is Selberg's integral:
$$
\int_{0}^{1} \dots \int_{0}^{1} \prod_{i=1}^{r} t_i^x (1-t_i)^y
\prod_{1 \leq i < j \leq r} (t_i-t_j)^{2z} dt_1 \dots dt_r
$$
$$
=\prod_{j=1}^{r}
{{(x+(j-1)z)!(y+(j-1)z)!(jz)!}
\over
{(x+y+(r+j-2)z+1)!z!}} \quad .
\eqno(Selberg)
$$
Now this is definitely a very nice identity. If you don't see why it
is nice, you did not yet develop the right aesthetic taste.
 
An even nicer identity is the contour-integral identity
$$
\int_{|z_1|=1} \dots \int_{|z_r|=1}  
\{\prod_{1 \leq i \neq j \leq r} (1- z_i/z_j)^a \} {{dz_1} \over {z_1}}
\dots {{dz_r} \over {z_r}} ={{(ra)!} \over {a!^r}} \quad .
\eqno(Dyson)
$$
[This identity has a good pedigree. It was conjectured
in 1960 by Dyson, and proved a year later by
Gusnson, and by Ken Wilson of Nobel (Renormalization Group) fame.
The `book' proof was given, in 1970, by the great
statistician I.J. Good.]
 
At present, using the multi-WZ method, Shalosh B. Ekhad can prove
$(Selberg)$,$(Dyson)$ (and any other such multi-dimensional identity),
for each specific number of variables $r$.
In practice, Shalosh can do it for $r=1,2,3,4$. It did.
Then the humans (in this case, myself and Wilf), looked at Shalosh's
proof, detected a common pattern, and it was a trivial,
albeit human, step, to formulate a WZ proof for
a general $r$ (see  our Inv. paper for the proof of $(Selberg)$).
 
BUT, right now, we still need this human factor. 
 
Another famous example is the Macdonald constant term
conjecture
(that Ian Macdonald talked about in his plenary talk,
in ICM 1998). Right now, Shalosh can do it for each
specific root system, but it takes a human of the
caliber of Cerednick to do it for {\it all root systems}.
 
But, NOT FOR LONG! Coming soon, in $\leq 10$ years is:
 
{\bf WZ, Chapter II (under construction)} 
 
In fact, {\it under construction} is wishful thinking.
So far I only have very few  epsilon-baked ideas. 
To wit:
 
Define precisely the notion of {\it hypergeometric
function} (what I called `nice') of $r$ variables. Where
$r$ is not {\it merely} a symbol {\it denoting} an integer,
but is a {\it symbol} period. Prototype nice functions
should be:
$$
\left ( \sum_{i=1}^{r} n_i \right )!
$$
$$
\prod_{i=1}^{r} NICE(n_i)
$$
and, probably
$$
\prod_{i=1}^{r-1} NICE(n_i+n_{i-1}) \quad ,
$$
$$
\prod_{1 \leq i < j \leq r} NICE(n_i+n_j) \quad ,
$$
where $NICE$ is a nice function of a single discrete variable.
 
The `language' should contain, as primitive symbols,
$\prod$ and $\sum$, that would have to be incorporated into
the algorithm.
 
We also need a notion of `global' niceness for
$NICE'_1 , \dots ,  NICE'_r$, i.e., it does not suffice
that $NICE'_i$ would be nice in its arguments,
but we should insist that $NICE_i'$, when
also viewed as a function of its {\it subscript}
$i$, is nice, in a sense yet to be made precise.
Somehow, we should also bring in symmetry.
 
So `soon', Macdonald's constant term identities, the Mehta integral,
the Selberg integral, and the Milne-Gustafson stuff would all be
fully automated.
 
{\bf Parable}
 
Prove the multi-variate identity
$$
\left ( \sum_{i=1}^{n} x_i \right )^2=
\sum_{i=1}^{n} x_i^2 + 2 \sum_{1 \leq i < j \leq n} x_ix_j \quad .
$$
 
{\bf Remark:} For each specific $n$, it is routine, but for
{\it general} $n$ we need a human, or do we?

{\bf Proof:} Let
$$
e_1(n):=\sum_{i=1}^{n} x_i \quad,
$$
$$
p_2(n):=\sum_{i=1}^{n} x_i^2 \quad,
$$
$$
e_2(n):=\sum_{1 \leq i < j \leq n} x_ix_j \quad.
$$
Now
$$
e_1(n)^2=(e_1(n-1)+x_n)^2=e_1(n-1)^2+2e_1(n-1)x_n+x_n^2
$$
$$
{ {induction} \atop {=} } p_2(n-1)+x_n^2+2e_2(n-1)+
2e_1(n-1)x_n= p_2(n) +2e_2(n) \quad.  \halmos 
$$
This is trivial to automate! We just have to work with
{\it indexed variables} $x[n], e1[n],e2[n],p2[n]$ etc.
 
Of course, the above identity is {\it not} WZ, but a
humdrum polynomial identity. BUT, 
one of the central messages of WZ theory is
that  `deep' hypergeometric summation and integration identities come down
to `trivial' rational function (and hence polynomial) identities.
So, it is very likely that the above parable has some substance.
\bye